\title{Graph powers, Delsarte, Hoffman, Ramsey and Shannon}
\author{{Noga Alon\thanks{School of Mathematics, Institute for Advanced Study, Princeton, NJ 08540, USA, and Raymond and Beverly Sackler
Faculty of Exact Sciences, Tel Aviv University, Tel Aviv, 69978,
Israel. Email: nogaa@tau.ac.il. Research supported in part by the
Israel Science Foundation, by a USA-Israeli BSF grant, by NSF grant
CCR-0324906, by a Wolfensohn fund and by the State of New Jersey.}}
\quad {Eyal Lubetzky
\thanks{ School of Computer Science, Raymond and Beverly
Sackler Faculty of Exact Sciences, Tel Aviv University, Tel Aviv,
69978, Israel. Email: lubetzky@tau.ac.il. Research partially
supported by a Charles Clore Foundation Fellowship.}}}
\documentclass [11pt]{article}
\usepackage[]{graphicx,epsf,epsfig,amsmath,amssymb,amsfonts,latexsym,amsthm}

\newtheorem{theorem}{Theorem}[section]
\newtheorem{lemma}[theorem]{Lemma}

\newtheorem*{definition}{Definition}

\newtheorem{corollary}[theorem]{Corollary}
\newtheorem{proposition}[theorem]{Proposition}

\renewcommand{\epsilon}{\varepsilon}

\newcommand{\xal}[1][p]{x_\alpha^{(#1)}}

\newtheoremstyle{upright}%
        {8pt plus2pt minus4pt}%
        {8pt plus2pt minus4pt}%
        {\upshape}%
        {}%
        {\bfseries}%
        {:}%
        {1em}%
        {}%
\theoremstyle{upright}
\newtheorem{remark}[theorem]{Remark}

\newcommand{\ignore}[1]{}

\begin{document}
\maketitle

\begin{abstract}
The $k$-th $p$-power of a graph $G$ is the graph on the vertex set
$V(G)^k$, where two $k$-tuples are adjacent iff the number of their
coordinates which are adjacent in $G$ is not congruent to $0$ modulo
$p$. The clique number of powers of $G$ is poly-logarithmic in the
number of vertices, thus graphs with small independence numbers in
their $p$-powers do not contain large homogenous subsets. We provide
algebraic upper bounds for the asymptotic behavior of independence
numbers of such powers, settling a conjecture of \cite{Xor} up to a
factor of $2$. For precise bounds on some graphs, we apply
Delsarte's linear programming bound and Hoffman's eigenvalue bound.
Finally, we show that for any nontrivial graph $G$, one can point
out specific induced subgraphs of large $p$-powers of $G$ with
neither a large clique nor a large independent set. We prove that
the larger the Shannon capacity of $\overline{G}$ is, the larger
these subgraphs are, and if $G$ is the complete graph, then some
$p$-power of $G$ matches the bounds of the Frankl-Wilson Ramsey
construction, and is in fact a subgraph of a variant of that
construction.
\end{abstract}

\section{Introduction}\label{sec::intro}
The $k$-th Xor graph power of a graph $G$, $G^{\oplus k}$, is the
graph whose vertex set is the cartesian product $V(G)^k$, where two
$k$-tuples are adjacent iff an odd number of their coordinates is
adjacent in $G$. This product was used in \cite{Thomason} to
construct edge colorings of the complete graph with two colors,
containing a smaller number of monochromatic copies of $K_4$ than
the expected number of such copies in a random coloring.

In \cite{Xor}, the authors studied the independence number,
$\alpha$, and the clique number, $\omega$, of high Xor powers of a
fixed graph $G$, motivated by problems in Coding Theory: cliques and
independent sets in such powers correspond to maximal codes
satisfying certain natural properties. It is shown in \cite{Xor}
that, while the clique number of $G^{\oplus k}$ is linear in $k$,
the independence number $\alpha(G^{\oplus k})$ grows exponentially:
the limit $\alpha(G^{\oplus k})^\frac{1}{k}$ exists, and is in the
range $[\sqrt{|V(G)|},|V(G)|]$. Denoting this limit by
$x_\alpha(G)$, the problem of determining $x_\alpha(G)$ for a given
graph $G$ proves to be extremely difficult, even for simple families
of graphs. Using spectral techniques, it is proved in \cite{Xor}
that $x_\alpha(K_n)=2$ for $n\in\{2,3,4\}$, where $K_n$ is the
complete graph on $n$ vertices, and it is conjectured that
$x_\alpha(K_n)=\sqrt{n}$ for every $n\geq 4$. The best upper bound
given in \cite{Xor} on $x_\alpha(K_n)$ for $n\geq 4$ is $n/2$.

The graph product we introduce in this work, which generalizes the
Xor product, is motivated by Ramsey Theory. In \cite{Erdos},
Erd\H{o}s proved the existence of graphs on $n$ vertices without
cliques or independent sets of size larger than $O(\log n)$
vertices, and that in fact, almost every graph satisfies this
property. Ever since, there have been many attempts to provide
explicit constructions of such graphs. Throughout the paper, without
being completely formal, we call a graph ``Ramsey" if it has neither
a ``large" clique nor a ``large" independent set. The famous Ramsey
construction of Frankl and Wilson \cite{FranklWilson} provided a
family of graphs on $n$ vertices, $FW_n$, with a bound of
$\exp\left(\sqrt{(2+o(1))\log n \log\log n}\right)$ on the
independence and clique numbers, using results from Extremal Finite
Set Theory. Thereafter, constructions with the same bound were
produced in \cite{NogaUnion} using polynomial spaces and in
\cite{Grolmusz} using low degree matrices. Recently, the old
Frankl-Wilson record was broken in \cite{BRSW}, where the authors
provide, for any $\epsilon>0$, a polynomial-time algorithm for
constructing a Ramsey graph on $n$ vertices without cliques or
independent sets on $\exp\left((\log n)^\epsilon\right)$ vertices.
The disadvantage of this latest revolutionary construction is that
it involves a complicated algorithm, from which it is hard to tell
the structure of the resulting graph.

Relating the above to graph products, the Xor product may be viewed
as an operator, $\oplus_k$, which takes a fixed input graph $G$ on
$n$ vertices, and produces a graph on $n^k$ vertices, $H=G^{\oplus
k}$. The results of \cite{Xor} imply that the output graph $H$
satisfies $\omega(H) \leq n k = O(\log(|V(H)|))$, and that if $G$ is
a nontrivial $d$-regular graph, then $H$ is $d'$-regular, with
$d'\to \frac{1}{2}|V(H)|$ as $k$ tends to infinity. Thus, $\oplus_k$
transforms any nontrivial $d$-regular graph into a random looking
graph, in the sense that it has an edge density of roughly
$\frac{1}{2}$ and a logarithmic clique number. However, the lower
bound $\alpha(H)\geq\sqrt{|V(H)|}$, which holds for every even $k$,
implies that $\oplus_k$ cannot be used to produce good Ramsey
graphs.

In order to modify the Xor product into a method for constructing
Ramsey graphs, one may try to reduce the high lower bound on the
independence numbers of Xor graph powers. Therefore, we consider a
generalization of the Xor graph product, which replaces the modulo
$2$ (adjacency of two $k$-tuples is determined by the parity of the
number of adjacent coordinates) with some possibly larger modulo
$p\in\mathbb{N}$. Indeed, we show that by selecting a larger $p$,
the lower bound on the independence number, $\alpha(H)$, is reduced
from $\sqrt{|V(H)|}$ to $|V(H)|^{1/p}$, at the cost of a polynomial
increase in $\omega(H)$. The generalized product is defined as
follows:
\begin{definition}
Let $k,p\in\mathbb{N}$. The $k$-th $p$-power of a graph $G$, denoted
by $G^{k_{(p)}}$, is the graph whose vertex set is the cartesian
product $V(G)^k$, where two $k$-tuples are adjacent iff the number
of their coordinates which are adjacent in $G$ is not congruent to
$0$ modulo $p$, that is:
$$(u_1,\ldots,u_k)~(v_1,\ldots,v_k)\in E(G^k) ~\mbox{ iff }~
|\{i:u_i v_i \in E(G)\}| \not\equiv 0\pmod{p}~.$$
\end{definition}

Throughout the paper, we use the abbreviation $G^k$ for
$G^{k_{(p)}}$ when there is no danger of confusion.

In Section \ref{sec::general-bounds} we show that the limit
$\alpha(G^k)^{\frac{1}{k}}$ exists and equals $\sup_k
\alpha(G^k)^{\frac{1}{k}}$; denote this limit by $\xal$. A simple
lower bound on $\xal$ is $|V(G)|^{1/p}$, and algebraic arguments
show that this bound is nearly tight for the complete graph:
$\xal(K_n)=O(n^{1/p})$. In particular, we obtain that $$\sqrt{n}
\leq x_\alpha(K_n) = \xal[2](K_n) \leq 2\sqrt{n-1}~,$$ improving the
upper bound of $n/2$ for $n\geq 4$ given in \cite{Xor}, and
determining that the behavior of $x_\alpha$ for complete graphs is
as stated in Question 4.1 of \cite{Xor} up to a factor of 2.

For the special case $G=K_n$, it is possible to apply Coding Theory
techniques in order to bound $\xal(G)$. The problem of determining
$\xal(K_n)$ can be translated into finding the asymptotic maximum
size of a code over the alphabet $[n]$, in which the Hamming
distance between any two codewords is divisible by $p$. The related
problem for {\em linear} codes over a field has been well studied:
see, e.g., \cite{WardDivisibleSurvey} for a survey on this subject.
However, as we later note in Section \ref{sec::general-bounds}, the
general non-linear case proves to be quite different, and the upper
bounds on linear divisible codes do not hold for $\xal(K_n)$. Yet,
other methods for bounding sizes of codes are applicable. In Section
\ref{sec::delsarte} we demonstrate the use of Delsarte's linear
programming bound in order to obtain precise values of
$\alpha(K_3^{k_{(3)}})$. We show that
$\alpha(K_3^{k_{(3)}})=3^{k/2}$ whenever $k\equiv 0\pmod{4}$, while
$\alpha(K_3^{k_{(3)}})< \frac{1}{2} 3^{k/2}$ for $k\equiv 2
\pmod{4}$, hence the series
$\alpha(K_3^{{k+1}_{(3)}})/\alpha(K_3^{k_{(3)}})$ does not converge
to a limit.

Section \ref{sec::hoffman} gives a general bound on $\xal$ for
$d$-regular graphs in terms of their eigenvalues, using Hoffman's
eigenvalue bound. The eigenvalues of $p$-powers of $G$ are
calculated using tensor products of matrices over $\mathbb{C}$, in a
way somewhat similar to performing a Fourier transform on the
adjacency matrix of $G$. This method may also be used to derive
tight results on $\alpha(G^{k_{(p)}})$, and we demonstrate this on
the above mentioned case of $p=3$ and the graph $K_3$, where we
compare the results with those obtained in Section
\ref{sec::delsarte} by the Delsarte bound.

Section \ref{sec::ramsey} shows, using tools from linear algebra,
that indeed the clique number of $G^{k_{(p)}}$ is poly-logarithmic
in $k$, and thus $p$-powers of graphs attaining the lower bound of
$\xal$ are Ramsey. We proceed to show a relation between the Shannon
capacity of the complement of $G$, $c(\overline{G})$, and the Ramsey
properties of $p$-powers of $G$. Indeed, for any nontrivial graph
$G$, we can point out a large Ramsey induced subgraph of some
$p$-power of $G$. The larger $c(\overline{G})$ is, the larger these
Ramsey subgraphs are. When $G=K_p$ for some prime $p$, we obtain
that $H=K_p^{{p^2}_{(p)}}$ is a Ramsey graph matching the bound of
Frankl-Wilson, and in fact, $H$ contains an induced subgraph which
is a modified variant of $FW_{N_1}$ for some $N_1$, and is contained
in another variant of $FW_{N_2}$ for some $N_2$. The method of
proving these bounds on $G^{k_{(p)}}$ provides yet another (simple)
proof for the Frankl-Wilson result.


\section{Algebraic lower and upper bounds on $x_\alpha^{(p)}$}\label{sec::general-bounds}
In this section, we define the parameter $\xal$, and provide lower
and upper bounds for it. The upper bounds follow from algebraic
arguments, using graph representation by polynomials.

\subsection{The limit of independence numbers of $p$-powers}
The following lemma establishes that $\xal$ exists, and gives simple
lower and upper bounds on its range for graphs on $n$ vertices:
\begin{lemma}\label{lem-xal-def}Let $G$ be a graph on $n$ vertices, and let
$p \geq 2$. The limit of $\alpha(G^{k_{(p)}})^{\frac{1}{k}}$ as
$k\to\infty$ exists, and, denoting it by $\xal(G)$, it satisfies:
$$ n^{1/p} \leq \xal(G) = \sup_k \alpha(G^{k_{(p)}})^{\frac{1}{k}} \leq n~.$$
\end{lemma}
\begin{proof}
Observe that if $I$ and $J$ are independent sets of $G^k$ and $G^l$
respectively, then the set $I \times J$ is an independent set of
$G^{k+l}$, as the number of adjacent coordinates between any two
$k$-tuples of $I$ and between any two $l$-tuples of $J$ is $0
\pmod{p}$. Therefore, the function $g(k)=\alpha(G^k)$ is
super-multiplicative and strictly positive, and we may apply
Fekete's Lemma (cf., e.g., \cite{VanLintWilson}, p. 85) to obtain
that the limit of $\alpha(G^k)^\frac{1}{k}$ as $k\to\infty$ exists,
and satisfies:
\begin{equation}\label{eq-xal-super-mult}\lim_{k\to\infty}
\alpha(G^k)^{\frac{1}{k}} = \sup_k
\alpha(G^k)^{\frac{1}{k}}~.\end{equation} Clearly, $\alpha(G^k) \leq
n^k$, and it remains to show the lower bound on $\xal$. Notice that
the following set is an independent set of $G^p$:
$$ I = \{~(u,\ldots,u)~:~u \in V(G)\} \subset G^p~,$$
since for all $u,v \in V(G)$, there are either $0$ or $p$ adjacent
coordinates between the two corresponding $p$-tuples in $I$. By
\eqref{eq-xal-super-mult}, we obtain that $\xal(G) \geq |I|^{1/p} =
n^{1/p}$.
\end{proof}

\subsection{Bounds on $\xal$ of complete graphs}
While the upper bound $|V(G)|$ on $\xal(G)$ is clearly attained by
an edgeless graph, proving that a family of graphs attains the lower
bound requires some effort. The next theorem states that complete
graphs achieve the lower bound of Lemma \ref{lem-xal-def} up to a
constant factor:
\begin{theorem}\label{thm-xal-kn-bounds}
The following holds for all integer $n,p \geq 2$:
\begin{equation}
  \label{eq-xal(kn)-upper-bound}
  \xal(K_n) \leq 2^{H(1/p)} (n-1)^{1/p}~,
\end{equation}
where $H(x)=-x\log_2(x)-(1-x)\log_2(1-x)$ is the binary entropy
function. In particular, $\xal(K_n)=\Theta(n^{1/p})$. In the special
case where $n=p=q^r$ for some prime $q$ and $r\geq 1$, the lower
bound roughly matches upper bound:
$$ p^{\frac{2}{p+1}} \leq \xal(K_p) \leq \left(\mathrm{e}p^2\right)^{1/p}~.$$
\end{theorem}
Taking $p=2$ and noting that $H(\frac{1}{2})=1$, we immediately
obtain the following corollary for Xor graph products, which
determines the asymptotic behavior of $x_\alpha$ for complete
graphs:
\begin{corollary}
For all $n \geq 2$, the complete graph on $n$ vertices satisfies
$$\sqrt{n} \leq x_\alpha(K_n) \leq 2\sqrt{n-1}~.$$
\end{corollary}
\begin{proof}
[Proof of Theorem \ref{thm-xal-kn-bounds}] The upper bound will
follow from an argument on polynomial representations, an approach
which was used in \cite{NogaUnion} to bound the Shannon capacity of
certain graphs. Take $k \geq 1$, and consider the graph $H = K_n^k$.
For every vertex of $H$, $u=(u_1,\ldots,u_k)$, we define the
following polynomial in $\mathbb{R}[x_{i,j}]$, where $i\in[k]$,
$j\in[n]$:
\begin{equation}
  \label{eq-poly-1-def-fu}
  f_u(x_{1,1},\ldots,x_{k,n}) = \prod_{t=1}^{\lfloor k/p \rfloor}
  \left(k-tp-\sum_{i=1}^k x_{i,u_i} \right)~.
\end{equation}
Next, give the following assignment of values for $\{x_{i,j}\}$,
$x_v$, to each $v=(v_1,\ldots,v_k)\in V(H)$:
\begin{equation}
  \label{eq-poly-1-def-x}
  x_{i,j} = \delta_{v_i,j}~,
\end{equation}
where $\delta$ is the Kronecker delta. Definitions
\eqref{eq-poly-1-def-fu} and \eqref{eq-poly-1-def-x} imply that for
every two such vertices $u=(u_1,\ldots,u_k)$ and
$v=(v_1,\ldots,v_k)$ in $V(H)$: \begin{equation}
\label{eq-poly-1-def-assignment}f_u(x_v) = \prod_{t=1}^{\lfloor k/p
\rfloor}
  \bigg(k-tp-\sum_{i=1}^k \delta_{u_i,v_i} \bigg)=
\prod_{t=1}^{\lfloor k/p \rfloor}
  \left(|\{i~:~u_i \neq v_i\}|-tp \right)
  ~.
\end{equation}
Notice that, by the last equation, $f_u(x_u)\neq 0$ for all $u \in
V(H)$, and consider two distinct non-adjacent vertices $u,v\in
V(H)$. The Hamming distance between $u$ and $v$ (considered as
vectors in $\mathbb{Z}_n^k$) is by definition $0 \pmod {p}$ (and is
not zero). Thus, \eqref{eq-poly-1-def-assignment} implies that
$f_u(x_v)=0$.

Recall that for all $u$, the assignment $x_u$ gives values
$x_{i,j}\in\{0,1\}$ for all $i,j$, and additionally, $\sum_{j=1}^n
x_{i,j}=1$ for all $i$. Therefore, it is possible to replace all
occurrences of $x_{i,n}$ by $1-\sum_{j=1}^{n-1} x_{i,j}$ in each
$f_u$, and then proceed and reduce the obtained result modulo the
polynomials:
$$ \bigcup_{i \in [k]} \left(\{ x_{i,j}^2-x_{i,j}:j\in[n]\}~ \cup ~\{ x_{i,j}
x_{i,l}:j,l\in [n], j \neq l\}\right)~,
$$ without affecting the value of the polynomials on the above defined
substitutions. In other words, after replacing $x_{i,n}$ by
$1-\sum_{j<n}x_{i,j}$, we repeatedly replace $x_{i,j}^2$ by
$x_{i,j}$, and let all the monomials containing $x_{i,j}x_{i,l}$ for
$j\neq l$ vanish. This gives a set of multi-linear polynomials
$\{\tilde{f}_u\}$ satisfying:
$$\left\{\begin{array}{cl}
\tilde{f}_u(x_u) \neq 0 &\mbox{for all }u \in V(H)\\
\tilde{f}_u(x_v) = 0 &\mbox{for }u\neq v~,~ u v \notin E(H)
\end{array}\right.~,$$
where the monomials of $\tilde{f}_u$ are of the form $\prod_{t=1}^r
x_{i_t,j_t}$ for some $0\leq r \leq \lfloor \frac{k}{p}\rfloor$, a
set of pairwise distinct indices $\{i_t\} \subset [k]$ and indices
$\{j_t\} \subset [n-1]$.

Let $\mathcal{F}= \mathrm{Span}(\{\tilde{f}_u:u\in V(H)\})$, and let
$I$ denote a maximum independent set of $H$. A standard argument
shows that $F=\{ \tilde{f}_u :u\in I\}$ is linearly independent in
$\mathcal{F}$. Indeed, suppose that $\sum_{u \in I} a_u
\tilde{f}_u(x) = 0$ ; then substituting $x=x_v$ for some $v\in I$
gives $a_v=0$. It follows that $\alpha(H) \leq \dim\mathcal{F}$, and
thus:
\begin{equation}\label{eq-poly-1-alpha-bound}
\alpha(H) \leq \sum_{r=0}^{\lfloor k/p \rfloor} \binom{k}{r} (n-1)^r
\leq  \left(2^{H(1/p)}(n-1)^{1/p}\right)^k ~,
\end{equation} where in the last inequality
we used the fact that $\sum_{i \leq \lambda n}\binom{n}{i} \leq 2^{n
H(\lambda)}$ (cf., e.g., the remark following Corollary 4.2 in
\cite{NogaExtSet}, and also \cite{ProbMethod} p. 242). Taking the
$k$-th root and letting $k$ grow to $\infty$, we obtain:
$$\xal(K_n) \leq 2^{H(1/p)}(n-1)^{1/p}~,$$ as required.

In the special case of $K_p$ (that is, $n=p$), note that:
$2^{H(\frac{1}{p})} = p^{\frac{1}{p}}(\frac{p}{p-1})^{\frac{p-1}{p}}
\leq (\mathrm{e}p)^{\frac{1}{p}}$ and hence in this case $\xal(K_p)
\leq (\mathrm{e}p^2)^{1/p}$. If $p=q^r$ is a prime-power we can
provide a nearly matching lower bound for $\xal(K_p)$ using a
construction of \cite{Xor}, which we shortly describe for the sake
of completeness.

Let $\mathcal{L}$ denote the set of all lines with finite slopes in
the affine plane $GF(p)$, and write down the following vector
$w_\ell$ for each $\ell \in \mathcal{L}$, $\ell =ax+b$ for some
$a,b\in GF(p)$:
$$ w_\ell = (a,a x_1 + b,a x_2+b,\ldots,a x_p+b)~,$$
where $x_1,\ldots,x_p$ denote the elements of $GF(p)$. For every two
distinct lines $\ell,\ell'$, if $\ell \| \ell'$ then
$w_{\ell},w_{\ell'}$ has a single common coordinate (the slope $a$).
Otherwise, $w_{\ell},w_{\ell'}$ has a single common coordinate,
which is the unique intersection of $\ell,\ell'$. In any case, we
obtain that the Hamming distance of $w_\ell$ and $w_{\ell'}$ is $p$,
hence $W=\{ w_\ell: \ell \in \mathcal{L}\}$ is an independent set in
$K_p^{p+1}$. By \eqref{eq-xal-super-mult}, we deduce that:
$$ \xal(K_p) \geq p^{\frac{2}{p+1}}~,$$
completing the proof.
\end{proof}
\begin{remark}
The proof of Theorem \ref{thm-xal-kn-bounds} used representation of
the vertices of $K_n^k$ by polynomials of $k n$ variables over
$\mathbb{R}$. It is possible to prove a similar upper bound on
$\xal(K_n)$ using a representation by polynomials of $k$ variables
over $\mathbb{R}$. To see this, use the natural assignment of $x_i =
v_i$ for $v=(v_1,\ldots,v_k)$, denoting it by $x_v$, and assign the
following polynomial to $u=(u_1,\ldots,u_k)$:
\begin{equation}\label{eq-poly-2-def-fu} f_u(x_1,\ldots,x_k)
= \prod_{t=1}^{\lfloor k/p \rfloor} \bigg( k-tp - \sum_{i=1}^k
\mathop{\prod_{j=1}^n}_{j \neq u_i}
\frac{x_i-j}{u_i-j}\bigg)~.\end{equation} The expression $\prod_{j
\neq u_i} \frac{x_i-j}{u_i-j}$ is the monomial of the Lagrange
interpolation polynomial, and is equal to $\delta_{x_i,u_i}$. Hence,
we obtain that $f_u(x_u)\neq 0$ for any vertex $u$, whereas
$f_u(x_v)=0$ for any two distinct non-adjacent vertices $u,v$. As
the Lagrange monomials yield values in $\{0,1\}$, we can convert
each $f_u$ to a multi-linear combination of these polynomials,
$\tilde{f}_u$, while retaining the above properties. Note that there
are $n$ possibilities for the Lagrange monomials (determined by the
value of $u_i$), and it is possible to express one as a linear
combination of the rest. From this point, a calculation similar to
that in Theorem \ref{thm-xal-kn-bounds} for the dimension of
$\mathrm{Span}(\{\tilde{f}_u:u\in V\})$ gives the upper bound
\eqref{eq-xal(kn)-upper-bound}.
\end{remark}

\begin{remark}
The value of $\alpha(K_n^{k_{(p)}})$ corresponds to a maximum size
of a code $C$ of $k$-letter words over $\mathbb{Z}_n$, where the
Hamming distance between any two codewords is divisible by $p$. The
case of {\em linear} such codes when $\mathbb{Z}_n$ is a field, that
is, we add the restriction that $C$ is a linear subspace of
$\mathbb{Z}_n^k$, has been thoroughly studied; it is equivalent to
finding a linear subspace of $\mathbb{Z}_n^k$ of maximal dimension,
such that the Hamming weight of each element is divisible by $p$. It
is known for this case that if $p$ and $n$ are relatively prime,
then the dimension of $C$ is at most $k/p$ (see
\cite{WardDivisible}), and hence the size of $C$ is at most
$n^{k/p}$. However, this bound does not hold for the non-linear case
(notice that this bound corresponds to the lower bound of Lemma
\ref{lem-xal-def}). We give two examples of this:
\begin{enumerate}
  \item Take $p=3$ and $n=4$. The divisible code bound implies an
  upper bound of $4^{1/3}\approx 1.587$, and yet $\xal[3](K_4) \geq
  \sqrt{3} \approx 1.732$. This follows from the geometric construction of Theorem
  \ref{thm-xal-kn-bounds}, which
  provides an independent set of size $9$ in $K_3^{4_{(3)}} \subset
  K_4^{4_{(3)}}$, using only the coordinates $\{0,1,2\}$ (this
  result can be slightly improved by adding an all-$3$ vector
  to the above construction in the $12$-th power).
  \item Take $p=3$ and $n=2$. The linear code bound is $2^{1/3}\approx
  1.26$, whereas the following construction shows that $\alpha(K_2^{12_{(3)}})\geq 24$, implying that $\xal[3](K_2) \geq
  {24}^{1/12}\approx 1.30$. Let $\{v_1,\ldots,v_{12}\}$ denote the rows
  of a binary Hadamard matrix of order 12 (such a matrix
  exists by Paley's Theorem, cf. e.g. \cite{MHall}).
  For all $i \neq j$, $v_i$ and $v_j$ have precisely $6$ common coordinates,
  and hence, the set
  $I=\{v_i\}\cup\{\overline{v}_i\}$ (where $\overline{v}_i$ denotes
  the complement of $v_i$ modulo 2) is an independent set of
  size 24 in $K_2^{12_{(3)}}$. In fact, $I$ is a maximum
  independent set of $K_2^{12_{(3)}}$, as Delsarte's linear programming
  bound (described in Section \ref{sec::delsarte}) implies that
  $\alpha(K_2^{12_{(3)}})\leq 24$.
\end{enumerate}

\end{remark}

\subsection{The value of $\xal[3](K_3)$}
While the upper bound of Theorem \ref{thm-xal-kn-bounds} on
$\xal(K_n)$ is tight up to a constant factor, the effect of this
constant on the independence numbers is exponential in the graph
power, and we must resort to other techniques in order to obtain
more accurate bounds. For instance, Theorem \ref{thm-xal-kn-bounds}
implies that:
$$ 1.732 \approx \sqrt{3} \leq \xal[3](K_3) \leq 2^{H(\frac{1}{3})}2^{\frac{1}{3}} =
\frac{3}{2^{1/3}} \approx 2.381~.
$$ In Sections \ref{sec::delsarte} and \ref{sec::hoffman}, we
demonstrate the use of Delsarte's linear programming bound and
Hoffman's eigenvalue bound for the above problem, and in both cases
obtain the exact value of $\alpha(K_3^{k_{(3)}})$ under certain
divisibility conditions. However, if we are merely interested in the
value of $\xal[3](K_3)$, a simpler consideration improves the bounds
of Theorem \ref{thm-xal-kn-bounds} and shows that
$\xal[3](K_3)=\sqrt{3}$:
\begin{lemma}
For any $k \geq 1$, $\alpha(K_3^{k_{(3)}}) \leq 3\cdot \sqrt{3}^k$,
and in particular, $\xal[3](K_3)=\sqrt{3}$.
\end{lemma}
\begin{proof}
Treating vertices of $K_3^k$ as vectors of $\mathbb{Z}_3^k$, notice
that every two vertices $x=(x_1,\ldots,x_k)$ and
$y=(y_1,\ldots,y_k)$ satisfy:
$$ \sum_{i=1}^k (x_i-y_i)^2 \equiv |\{ i: x_i \neq y_i\}| \pmod{3}~,$$
and hence if $I$ is an independent set in $K_3^k$, then:
$$ \sum_i (x_i-y_i)^2 \equiv 0 \pmod{3} ~\mbox{ for all }
x,y \in I~.$$ Let $I$ denote a maximum independent set of $K_3^k$,
and let $I_c = \{ x\in I: \sum_i x_i^2 \equiv c \pmod{3}\}$ for
$c\in\{0,1,2\}$. For every $c\in\{0,1,2\}$ we have:
$$ \sum_i (x_i-y_i)^2 = 2c - 2 x\cdot y \equiv 0
\pmod{3}~\mbox{ for all }x,y\in I_c,$$ and hence $x \cdot y = c$ for
all $x,y\in I_c$. Choose $c$ for which $|I_c| \geq |I|/3$, and
subtract an arbitrary element $z \in I_c$ from all the elements of
$I_c$. This gives a set $J$ of size at least $|I|/3$, which
satisfies:
$$ x \cdot y = 0 ~\mbox{ for all }x,y\in J~.$$ Since
$\mathrm{Span}(J)$ is a self orthogonal subspace of
$\mathbb{Z}_3^k$, its dimension is at most $k/2$, and hence $|J|\leq
3^{k/2}$. Altogether, $\alpha(K_3^k) \leq 3 \cdot \sqrt{3}^k$, as
required.
\end{proof}

\section{Delsarte's linear programming bound for complete graphs}\label{sec::delsarte}
In this section, we demonstrate how Delsarte's linear programming
bound may be used to derive precise values of independence numbers
in $p$-powers of complete graphs. As this method was primarily used
on binary codes, we include a short proof of the bound for a general
alphabet.

\subsection{Delsarte's linear programming bound}

The linear programming bound follows from the relation between the
distance distribution of codes and the Krawtchouk polynomials,
defined as follows:
\begin{definition}
Let $n\in\mathbb{N}$ and take $q \geq 2$. The Krawtchouk polynomials
$\mathcal{K}_k^{n;q}(x)$ for $k=0,\ldots,n$ are defined by:
\begin{equation}\label{eq-kraw-def} \mathcal{K}_k^{n;q}(x) = \sum_{j=0}^k \binom{x}{j}
\binom{n-x}{k-j} (-1)^j (q-1)^{k-j}~.\end{equation}
\end{definition}
\begin{definition}
Let $C$ be an $n$-letter code over the alphabet $\{1,\ldots,q\}$.
The distance distribution of $C$, $B_0,B_1,\ldots,B_n$, is defined
by:
$$ B_k = \frac{1}{|C|}|\{ (w_1,w_2)\in C^2 : \delta(w_1,w_2)=
k\}|~~(k=0,\ldots,n)~,$$ where $\delta$ denotes the Hamming
distance.
\end{definition}
The Krawtchouk polynomials $\{\mathcal{K}_k^{n;q}(x)\}$ are
sometimes defined with a normalizing factor of $q^{-k}$. Also, it is
sometimes customary to define the distance distribution with a
different normalizing factor, letting $A_k = \frac{B_k}{|C|}$, in
which case $A_k$ is the probability that a random pair of codewords
has a Hamming distance $k$.

The Krawtchouk polynomials $\{ \mathcal{K}_k^{n;q} : k=0,\ldots,n\}$
form a system of orthogonal polynomials with respect to the weight
function $w(x)= \frac{n!}{\Gamma(1+x)\Gamma(n+1-x)}(q-1)^x$, where
$\Gamma$ is the gamma function. For further information on these
polynomials see, e.g., \cite{Szego}.

Delsarte \cite{Delsarte1} (see also \cite{MS}) presented a
remarkable method for bounding the maximal size of a code with a
given set of restrictions on its distance distribution. This
relation is given in the next proposition, for which we include a
short proof:
\begin{proposition}\label{prop-delsarte-relation}
Let $C$ be a code of $n$-letter words over the alphabet $[q]$, whose
distance distribution is $B_0,\ldots,B_n$. The following holds:
\begin{equation}\label{eq-kraw-distances-relation} \sum_{i=0}^n B_i
\mathcal{K}_k^{n;q}(i) \geq 0 ~\mbox{ for all
}k=0,\ldots,n~.\end{equation}
\end{proposition}
\begin{proof}
Let $G=\mathbb{Z}_q^n$, and for every two functions $f,g:G\to
\mathbb{C}$, define (as usual) their inner product
$\left<f,g\right>$ and their delta-convolution, $f*g$, as:
$$ \left<f,g\right> = \int_G f(x) \overline{g(x)} dx = \frac{1}{|G|}\sum_{T \in G}f(T)
\overline{g(T)}~,$$
$$(f * g )(s) = \int_G f(x) \overline{g(x-s)} dx~.$$
Denoting the Fourier expansion of $f$ by: $f=\sum_{S\in G}
\widehat{f}(S) \chi_S$, where $\chi_S(x)=\omega^{S \cdot x}$ and
$\omega$ is the $q$-th root of unity, it follows that for any
$k=0,\ldots,n$:
\begin{equation}
  \label{eq-kraw-fourier-formula}
  \sum_{S\in G : |S|=k}\widehat{f}(S) = \frac{1}{|G|}\sum_{i=0}^n
  \mathcal{K}_k^{n;q}(i) \sum_{T \in G : |T|=i} f(T)~,
\end{equation}
where $|S|$ and $|T|$ denote the Hamming weights of $S,T \in G$.
Since the delta-convolution satisfies: $$\widehat{f*g}(S) =
\widehat{f}(S)\overline{\widehat{g}(S)}~,$$ every $f$ satisfies:
\begin{equation}
  \label{eq-characteristic-convolution-geq-0}
  \widehat{f*f}(S) = |\widehat{f}(S)|^2 \geq 0~.
\end{equation} Let $f$ denote the
characteristic function of the code $C$, $f(x)=\mathbf{1}_{\{x\in
C\}}$, and notice that:
$$ (f*f)(S) = \frac{1}{|G|}\sum_{T\in G}f(T)\overline{f(T-S)} = \frac{1}{|G|}|\{T: T,T-S\in C\}|~,$$
and thus:
\begin{equation}
  \label{eq-bi-convolution-relation}
B_i = \frac{|G|}{|C|}\sum_{T:|T|=i}(f*f)(T)~.
\end{equation}
Putting together \eqref{eq-kraw-fourier-formula},
\eqref{eq-characteristic-convolution-geq-0} and
\eqref{eq-bi-convolution-relation}, we obtain: $$ 0 \leq
\sum_{S:|S|=k}\widehat{f*f}(S) = \frac{1}{|G|} \sum_{i=0}^n
\mathcal{K}_k^{n;q}(i) \sum_{T:|T|=i}(f*f)(T) =
\frac{|C|}{|G|^2}\sum_{i=0}^n \mathcal{K}_k^{n;q}(i) B_i ~,$$ as
required.
\end{proof}

Let $F\subset [n]$ be a set of forbidden distances between distinct
codewords. Since $|C|=\sum_i B_i$, the following linear program
provides an upper bound on the size of any code with no pairwise
distances specified by $F$:
$$\begin{array}{c}
\mbox{maximize }\sum_i B_i \mbox{ subject to the constraints: }\\
\left\{\begin{array}{c}B_0=1 \\ B_i \geq 0 \mbox{ for
all }i\\ B_i = 0 \mbox{ for all } i \in F \\
\sum_{i=0}^n B_i \mathcal{K}_k^{n;q}(i) \geq 0 ~\mbox{ for all
}k=0,\ldots,n\end{array}\right.\end{array}~.$$ By examining the dual
program, it is possible to formulate this bound as a minimization
problem. The following proposition has been proved in various
special cases, (cf., e.g., \cite{Delsarte2}, \cite{Litsyn}). For the
sake of completeness, we include a short proof of it.
\begin{proposition}\label{prop-delsarte-bound} Let $C$ be a code of $n$-letter words over the
alphabet $[q]$, whose distance distribution is $B_0,\ldots,B_n$. Let
$\displaystyle{P(x)=\sum_{k=0}^n \alpha_k \mathcal{K}_k^{n;q}(x)}$
denote an $n$-degree polynomial over $\mathbb{R}$. If $P(x)$ has the
following two properties:
\begin{equation}\label{eq-delsarte-poly-req-1}
\alpha_0 > 0 ~\mbox{ and }\alpha_i \geq 0 ~\mbox{ for all }i
=1,\ldots,n~,
\end{equation}
\begin{equation}\label{eq-delsarte-poly-req-2}
P(d) \leq 0 ~\mbox{ whenever }B_d > 0 \mbox { for }d=1,\ldots,n~,
\end{equation}
then $|C| \leq P(0)/\alpha_0$.
\end{proposition}
\begin{proof}
The Macwilliams transform of the vector $(B_0,\ldots,B_n)$ is
defined as follows:
\begin{equation}\label{eq-bk'-definition}
B'_k = \frac{1}{|C|}\sum_{i=0}^n
\mathcal{K}_k^{n;q}(i)B_i~.\end{equation} By the Delsarte
inequalities (stated in Proposition \ref{prop-delsarte-relation}),
$B'_k \geq 0$, and furthermore:
$$B'_0 = \frac{1}{|C|}\sum_{i=0}^n \mathcal{K}_0^{n;q}(i) B_i =
\frac{1}{|C|}\sum_i B_i = 1~.$$ Therefore, as
\eqref{eq-delsarte-poly-req-1} guarantees that $\alpha_i \geq 0$ for
$i> 0$, we get:
\begin{equation} \label{eq-bk'-lower-bound} \sum_{k=0}^n \alpha_k
B'_k \geq \alpha_0~.
\end{equation}
On the other hand, $B_0=1$, and by \eqref{eq-delsarte-poly-req-2},
whenever $B_i > 0$ for some $i>0$ we have $P(i)\leq 0$, thus:
\begin{equation}
  \label{eq-bi-Pi-upper-bound} \sum_{i=0}^n B_i P(i) \leq P(0)~.
\end{equation}
Combining \eqref{eq-bk'-lower-bound} and
\eqref{eq-bi-Pi-upper-bound} with \eqref{eq-bk'-definition} gives:
$$
\alpha_0 \leq \sum_{k=0}^n \alpha_k B'_k = \frac{1}{|C|}
\sum_{i=0}^n B_i \sum_{k=0}^n \alpha_k \mathcal{K}_k^{n;q}(i) =
\frac{1}{|C|} \sum_{i=0}^n B_i P(i) \leq \frac{P(0)}{|C|}~,
$$
and the result follows.
\end{proof}

We proceed with an application of the last proposition in order to
bound the independence numbers of $p$-powers of complete graphs. In
this case, the distance distribution is supported by $\{ i : i
\equiv 0 \pmod{p}\}$, and in Section \ref{sec::delsarte-application}
we present polynomials which satisfy the properties of Proposition
\ref{prop-delsarte-bound} and provide tight bounds on
$\alpha(K_3^{k_{(3)}})$.

\subsection{Improved estimations of $\alpha(K_3^{k_{(3)}})$}\label{sec::delsarte-application}
Recall that the geometric construction of Theorem
\ref{thm-xal-kn-bounds} describes an independent set of size $p^2$
in $K_p^{{p+1}_{(p)}}$ for every $p$ which is a prime-power. In
particular, this gives an independent set of size $3^{k/2}$ in
$K_3^{k_{(3)}}$ for every $k\equiv 0\pmod{4}$. Using Proposition
\ref{prop-delsarte-bound} we are able to deduce that indeed
$\alpha(K_3^k)=3^{k/2}$ whenever $k\equiv 0\pmod{4}$, whereas for
$k\equiv 2\pmod{4}$ we prove that $\alpha(K_3^k)< \frac{1}{2}
3^{k/2}$.
\begin{theorem}\label{thm-delsarte-k3-bounds}
The following holds for any even integer $k$:
$$ \left\{\begin{array}{ll}
\alpha(K_3^k) = 3^{k/2} & k \equiv 0 \pmod{4} \\
\frac{1}{3} 3^{k/2} \leq \alpha(K_3^k)  < \frac{1}{2} 3^{k/2} & k
\equiv 2 \pmod{4}
\end{array}\right. ~.$$
\end{theorem}
\begin{proof}
Let $k$ be an even integer, and define the following polynomials:
\begin{eqnarray}
  P(x) = \frac{2}{3} 3^{k/2} + \mathop{\sum_{t=1}^k}_{t \not\equiv 0
  (\mathrm{mod}~3)} \mathcal{K}_t^{k;3}(x)~, \label{eq-delsarte-poly-1} \\
  Q(x) = \frac{2}{3} 3^{k/2} + \mathop{\sum_{t=0}^k}_{t \equiv 0
  (\mathrm{mod}~3)} \mathcal{K}_t^{k;3}(x)~. \label{eq-delsarte-poly-2}
\end{eqnarray}
Clearly, both $P$ and $Q$ satisfy \eqref{eq-delsarte-poly-req-1}, as
$\mathcal{K}_0^{n;q}=1$ for all $n,q$. It remains to show that $P,Q$
satisfy \eqref{eq-delsarte-poly-req-2} and to calculate $P(0),Q(0)$.
As the following calculation will prove useful later on, we perform
it for a general alphabet $q$ and a general modulo $p$. Denoting the
$q$-th root of unity by $\omega=\mathrm{e}^{2\pi i/q}$, we have:
\begin{eqnarray}
\mathop{\sum_{t=0}^k}_{t \equiv 0 (\mathrm{mod}~p)}
\mathcal{K}_t^{k;q}(s) &=& \mathop{\sum_{t=0}^k}_{t \equiv 0
(\mathrm{mod}~p)} \sum_{j=0}^t
\binom{s}{j} \binom{k-s}{t-j}(-1)^j (q-1)^{t-j}= \nonumber \\
&=& \sum_{j=0}^s \binom {s}{j} (-1)^j
\mathop{\sum_{l=0}^{k-s}}_{j+l\equiv 0(\mathrm{mod}~p)}
\binom{k-s}{l} (q-1)^l = \nonumber\\
 &=& \sum_{j=0}^s \binom {s}{j}
(-1)^j \sum_{l=0}^{k-s} \binom{k-s}{l} (q-1)^l \frac{1}{q}
\sum_{t=0}^{q-1}
\omega^{(j+l)t}= \nonumber \\
&=& \delta_{s,0} \cdot q^{k-1} + \frac{1}{q}\sum_{t=1}^{q-1}
(1+(q-1)\omega^t)^{k-s} (1-\omega^t)^s~, \label{eq-Kk-0-mod-p-sum}
\end{eqnarray} where the last equality is by the fact that:
$\sum_{j=0}^s\binom{s}{j}(-1)^j = \delta_{s,0}$, and therefore the
summand for $t=0$ vanishes if $s\neq 0$ and is equal to $q^{k-1}$ if
$s=0$. Repeating the above calculation for $t\not\equiv 0
\pmod{p}$ gives:\begin{eqnarray} \mathop{\sum_{t=0}^k}_{t \not\equiv
0 (\mathrm{mod}~p)} \mathcal{K}_t^{k;q}(s) &=& \sum_{j=0}^s \binom
{s}{j} (-1)^j \sum_{l=0}^{k-s} \binom{k-s}{l} (q-1)^l
\left(1-\frac{1}{q} \sum_{t=0}^{q-1} \omega^{(j+l)t}\right)=
\nonumber \\
&=& \delta_{s,0}\cdot (q^k - q^{k-1}) -\frac{1}{q}\sum_{t=1}^{q-1}
(1+(q-1)\omega^t)^{k-s}
(1-\omega^t)^s\label{eq-Kk-non-0-mod-p-sum}~.
\end{eqnarray}
Define:
$$ \xi_s = \frac{1}{q}\sum_{t=1}^{q-1}
(1+(q-1)\omega^t)^{k-s} (1-\omega^t)^s ~,$$ and consider the special
case $p=q=3$. The fact that $\omega^2=\overline{\omega}$ implies
that: \begin{equation}\label{eq-xi-s-value} \xi_s =
\frac{2}{3}\mathrm{Re}\left((1+2\omega)^{k-s}(1-\omega)^s\right)=
\frac{2}{3}\mathrm{Re}\left((\sqrt{3}i)^{k-s}
(\sqrt{3}\mathrm{e}^{-\frac{\pi}{6}i})^s\right)
=\frac{2}{3}\sqrt{3}^k\cos(\frac{\pi k}{2}-\frac{2\pi s}{3})~,
\end{equation} and for even values of $k$ and $s \equiv 0\pmod{3}$ we deduce
that:
\begin{equation}\label{eq-Kk-exp-p-q-3}\xi_s = \frac{2}{3} 3^{k/2} (-1)^{k/2}~.\end{equation}
 Therefore, $\xi_s = \frac{2}{3}3^{k/2}$ whenever $s\equiv
0\pmod{3}$ and $k\equiv 0\pmod{4}$, and
\eqref{eq-Kk-non-0-mod-p-sum} gives the following for any $k \equiv
0 \pmod{4}$:
\begin{eqnarray} P(0) &=& \frac{2}{3} 3^{k/2} + \frac{2}{3} 3^k -
\xi_0 = \frac{2}{3}
3^k~,\nonumber\\
P(s) &=& \frac{2}{3} 3^{k/2} -\xi_s = 0~\mbox{ for any }0\neq s
\equiv 0\pmod{3}~.\nonumber\end{eqnarray} Hence, $P(x)$ satisfies
the requirements of Proposition \ref{prop-delsarte-bound} and we
deduce that for any $k\equiv 0\pmod{4}$:
$$ \alpha(K_3^k) \leq \frac{P(0)}{\frac{2}{3}3^{k/2}} = 3^{k/2}~.$$
As mentioned before, the construction used for the
lower bound on $\xal(K_3)$ implies that this bound is indeed tight
whenever $4 \mid k$.

For $k \equiv 2 \pmod{4}$ and $s\equiv 0 \pmod{3}$ we get $\xi_s
=-\frac{2}{3}3^{k/2}$, and by \eqref{eq-Kk-0-mod-p-sum} we get:
\begin{eqnarray} Q(0) &=& \frac{2}{3} 3^{k/2} + 3^{k-1} +
\xi_0 = 3^{k-1}~,\nonumber\\
Q(s) &=& \frac{2}{3} 3^{k/2} + \xi_s = 0~\mbox{ for any }0\neq s
\equiv 0\pmod{3}~.\nonumber\end{eqnarray} Again, $Q(x)$ satisfies
the requirements of Proposition \ref{prop-delsarte-bound} and we
obtain the following bound for $k\equiv 2\pmod{4}$:
$$ \alpha(K_3^k) \leq \frac{Q(0)}{
\frac{2}{3}3^{k/2}+1} = \frac{3^k}{2 \cdot 3^{k/2} + 3} <
\frac{1}{2} 3^{k/2}~.$$ To conclude the proof, take a maximum
independent set of size $\sqrt{3}^l$ in $K_3^l$, where $l=k-2$, for
a lower bound of $\frac{1}{3}3^{k/2}$.

\end{proof}

\section{Hoffman's bound on independence numbers of $p$-powers}\label{sec::hoffman}
In this section we apply spectral analysis in order to bound the
independence numbers of $p$-powers of $d$-regular graphs. The next
theorem generalizes Theorem 2.9 of \cite{Xor} by considering tensor
powers of adjacency matrices whose values are $p$-th roots of unity.
\begin{theorem}\label{thm-hoffman-bound}Let $G$ be a nontrivial $d$-regular graph on $n$
vertices, whose eigenvalues are $d=\lambda_1\geq \lambda_2 \geq
\ldots \geq \lambda_n$, and let
$\lambda=\max\{\lambda_2,|\lambda_n|\}$. The following holds for any
$p\geq 2$:
\begin{equation}\label{eq-hoff-xal-bound}
\xal(G) \leq \max \{
\sqrt{n^2-2\left(1-\cos(\frac{2\pi}{p})\right)d(n-d)}, \lambda
\sqrt{2-2\cos\left(\frac{2\pi}{p}\lfloor\frac{p}{2}\rfloor\right)}
\}~.
\end{equation}
\end{theorem}
\begin{proof}
Let $A = A_G$ denote the adjacency matrix of $G$, and define the
matrices $B_t$ for $t\in \mathbb{Z}_p$ as follows:
\begin{equation}\label{eq-B-t-definition}
B_t = J_n + (\omega^t -1)A~,
\end{equation}
where $\omega=\mathrm{e}^{2 \pi i/p}$ is the $p$-th root of unity,
and $J_n$ is the all-ones matrix of order $n$. In other words:
$$ (B_t)_{u v} = \omega^{t A_{u v}} = \left\{\begin{array}
  {ll}\omega^t & \mbox{if }u v \in E(G) \\
  1 &\mbox{if } u v \notin E(G)
\end{array}\right.~.$$
By the definition of the matrix tensor product $\otimes$, it follows
that for all $u=(u_1,\ldots,u_k)$ and $v=(v_1,\ldots,v_k)$ in $G^k$:
$$ (B_t^{\otimes k})_{u,v} = \omega^{t
 |\{i ~:~ u_i v_i \in E(G)\}|}~,$$
 and:
$$ \sum_{t=0}^{p-1} (B_t^{\otimes k})_{u,v} = \left\{\begin{array}
  {ll}p & \mbox{if } ~|\{i : u_i v_i \in E(G)\}| \equiv 0 \pmod{p}\\
  0 &\mbox{otherwise}
\end{array}\right.~.$$
Recalling that $u v \in E(G^k)$ iff $|\{i: u_i v_i \in
E(G)\}|\not\equiv 0 \pmod{p}$, we get:
\begin{equation}
  \label{eq-A-G^k-B-relation}
  A_{G^k} = J_{n^k} - \frac{1}{p}\sum_{t=0}^{p-1}B_t^{\otimes k} =
  \frac{p-1}{p} J_{n^k} - \frac{1}{p}\sum_{t=1}^{p-1}B_t^{\otimes k}~.
\end{equation}
The above relation enables us to obtain expressions for the
eigenvalues of $G^k$, and then apply the following bound, proved by
Hoffman (see \cite{Hoffman}, \cite{LovaszTheta}): every regular
nontrivial graph $H$ on $N$ vertices, whose eigenvalues are $\mu_1
\geq \ldots \geq \mu_N$, satisfies:
\begin{equation}\label{eq-hoffman-bound}
  \alpha(H) \leq \frac{-N \mu_N}{\mu_1-\mu_N}~.
\end{equation}
Recall that $J_n$ has a single non-zero eigenvalue of $n$,
corresponding to the all-ones vector $\underline{1}$. Hence,
\eqref{eq-B-t-definition} implies that $\underline{1}$ is an
eigenvector of $B_t$ with an eigenvalue of $n+(\omega^t-1)d$, and
the remaining eigenvalues of $B_t$ are
$\{(\omega^t-1)\lambda_i:i>1\}$. By well known properties of tensor
products, we obtain that the largest eigenvalue of $H=G^k$ (which is
its degree of regularity) is:
\begin{eqnarray} \mu_1 &=& n^k -
\frac{1}{p}\sum_{t=0}^{p-1}(n+(\omega^t-1)d)^k =
 n^k - \frac{1}{p}\sum_{j=0}^k \binom{k}{j}(n-d)^{k-j} d^j \sum_{t=0}^{p-1}\omega^{j
 t}
 = \nonumber \\
 &=& n^k - \mathop{\sum_{j=0}^k}_{j\equiv 0
(\mathrm{mod}~p)} \binom{k}{j}(n-d)^{k-j} d^j ~,\label{eq-mu1-d-reg}
\end{eqnarray}
and the remaining eigenvalues are of the form:
\begin{equation}\label{eq-mu}
\mu(\lambda_{i_1},\ldots,\lambda_{i_s}) = -\frac{1}{p}
\sum_{t=1}^{p-1}(n+(\omega^t-1)d)^{k-s}\prod_{j=1}^s
(\omega^t-1)\lambda_{i_j}~,\end{equation} where $0<s\leq k$ and $1 <
i_j \leq n$ for all $j$ (corresponding to an eigenvector which is a
tensor-product of the eigenvectors of $\lambda_{i_j}$ for
$j=1,\ldots,s$ and $\underline{1}^{\otimes k-s}$). The following
holds for all such choices of $s$ and $\{\lambda_{i_j}\}$:
\begin{eqnarray} |\mu (\lambda_{i_1},\ldots,\lambda_{i_s})| &\leq& \max_{1\leq t \leq
p-1} \bigg|(n+(\omega^t-1)d)^{k-s}
\prod_{i=1}^s(\omega^t-1)\lambda_{i_j} \bigg| \leq
\nonumber \\
& \leq & \max_{1\leq t \leq p-1} |n+(\omega^t-1)d|^{k-s}
(|\omega^t-1|\lambda)^s \leq \nonumber \\
& \leq & \max_{1\leq t \leq p-1}
\left(\max\{|n+(\omega^t-1)d|,\lambda|\omega^t-1|\}\right)^k~.
\nonumber
\end{eqnarray}
Since for any $1 \leq t \leq p-1$ we have:
\begin{eqnarray} &|n+(\omega^t-1)d|^2 &= n^2 - 2\left(1-\cos(\frac{2\pi
t}{p})\right)d(n-d) \leq
n^2-2\left(1-\cos(\frac{2\pi}{p})\right)d(n-d) ~,\nonumber \\
&|\omega^t-1|^2 &= 2-2\cos(\frac{2\pi t}{p}) \leq
2-2\cos\left(\frac{2\pi}{p}\lfloor\frac{p}{2}\rfloor\right) ~,
\nonumber \end{eqnarray} it follows that:
$$ |\mu(\lambda_{i_1},\ldots,\lambda_{i_s})| \leq ( \max\{ \rho_1, \rho_2 \} )^k~,$$
where: $$\begin{array}{lll}
\rho_1 &=&\sqrt{n^2-2\left(1-\cos(\frac{2\pi}{p})\right)d(n-d)}\\
\rho_2
&=&\lambda\sqrt{2-2\cos\left(\frac{2\pi}{p}\lfloor\frac{p}{2}\rfloor\right)}\end{array}~.$$
By the same argument, \eqref{eq-mu1-d-reg} gives:
$$ |\mu_1| \geq n^k - \rho_1^k~,$$
and applying Hoffman's bound \eqref{eq-hoffman-bound}, we get:
\begin{equation}\label{eq-hoffman-rho-bound}\alpha(G^k) \leq \frac{-n^k \mu_{n^k}}{\mu_1 -
\mu_{n^k}} \leq \frac{(\max\{\rho_1,\rho_2\})^k }{1 -
(\frac{\rho_1}{n})^k +
(\frac{\max\{\rho_1,\rho_2\}}{n})^k}~.\end{equation} To complete the
proof, we claim that $\max\{\rho_1,\rho_2\} \leq n$, and hence the
denominator in the expression above is $\Theta(1)$ as $k\to\infty$.
Clearly, $\rho_1 \leq n$, and a simple argument shows that $\lambda
\leq n/2$ and hence $\rho_2 \leq n$ as well. To see this, consider
the matrix $A^2$ whose diagonal entries are $d$; we have:
$$ n d =
\mathrm{tr} A^2 = \sum_i \lambda_i^2 \geq d^2 + \lambda^2 ~,$$
implying that $\lambda \leq \sqrt{d(n-d)} \leq \frac{n}{2}$.
Altogether, taking the $k$-th root and letting $k$ tend to $\infty$
in \eqref{eq-hoffman-rho-bound}, we obtain that $\xal(G) \leq
\max\{\rho_1,\rho_2\}$, as required.
\end{proof}
\noindent \textbf{Examples:} For $p=2,3$ the above theorem gives:
\begin{eqnarray} \xal[2](G) &\leq&
\max\{|n-2d|,2\lambda\}~,\nonumber\\
\xal[3](G) &\leq&
\max\{\sqrt{n^2-3d(n-d)},\sqrt{3}\lambda\}~.\nonumber\end{eqnarray}

Since the eigenvalues of $K_3$ are $\{2,-1,-1\}$, this immediately
provides another proof for the fact that $\xal[3](K_3) \leq
\sqrt{3}$. Note that, in general, the upper bounds derived in this
method for $\xal(K_n)$ are only useful for small values of $n$, and
tend to $n$ as $n\to\infty$, whereas by the results of Section
\ref{sec::general-bounds} we know that $\xal(K_n)=\Theta(n^{1/p})$.

Consider $d=d(n)=\frac{n}{2}+O(\sqrt{n})$, and let $G \sim G_{n,d}$
denote a random $d$-regular graph on $n$ vertices. By the results of
\cite{KSVW}, $\lambda = \max\{\lambda_2,|\lambda_n|\} = O(n^{3/4})$,
and thus, Theorem \ref{thm-hoffman-bound} implies that $\xal[2](G) =
O(n^{3/4})$, and $\xal[3](G) \leq (1+o(1))\frac{n}{2}$. We note that
one cannot hope for better bounds on $\xal[3]$ in this method, as
$\rho_1$ attains its minimum at $d=\frac{n}{2}$.
\begin{remark}
The upper bound \eqref{eq-hoff-xal-bound} becomes weaker as $p$
increases. However, if $p$ is divisible by some $q \geq 2$, then
clearly any independent set of $G^{k_{(p)}}$ is also an independent
set of $G^{k_{(q)}}$, and in particular, $\xal[p](G) \leq
\xal[q](G)$. Therefore, when applying Theorem
\ref{thm-hoffman-bound} on some graph $G$, we can replace $p$ by the
minimal $q \geq 2$ which divides $p$. For instance, $\xal[4](G) \leq
\xal[2](G)\leq \max\{|n-2d|,2\lambda\}$, whereas substituting $p=4$
in \eqref{eq-hoff-xal-bound} gives the slightly weaker bound
$\xal[4](G)\leq \{\sqrt{(n-d)^2+d^2},2\lambda\}$.
\end{remark}
\begin{remark}
In the special case $G=K_n$, the eigenvalues of $G$ are
$\{n-1,-1,\ldots,-1\}$, and the general expression for the
eigenvalues of $G^k$ in \eqref{eq-mu} takes the following form (note
that $\lambda_{i_j}=-1$ for all $1\leq j \leq s$):
$$ \mu(s) = -\frac{1}{p}
\sum_{t=1}^{p-1}(1+(n-1)\omega^t)^{k-s} (1-\omega^t)^s~,$$ and as
$s>0$, we obtain the following from \eqref{eq-Kk-non-0-mod-p-sum}:
$$ \mu(s) =  \mathop{\sum_{t=0}^k}_{t \not\equiv 0 (\mathrm{mod}~p)}
\mathcal{K}_t^{k;q}(s) ~.$$ Similarly, comparing
\eqref{eq-mu1-d-reg} to \eqref{eq-Kk-non-0-mod-p-sum} gives:
$$ \mu_1 =  \mathop{\sum_{t=0}^k}_{t \not\equiv 0 (\mathrm{mod}~p)}
\mathcal{K}_t^{k;q}(0) ~.$$ It is possible to deduce this result
directly, as $K_n^k$ is a Cayley graph over $\mathbb{Z}_n^k$ with
the generator set $ S = \{x : |x| \not\equiv 0 \pmod{p}\}$, where
$|x|$ denotes the Hamming weight of $x$. It is well known that the
eigenvalues of a Cayley graph are equal to the character sums of the
corresponding group elements. Since for any $k=0,\ldots,n$ and any
$x \in \mathbb{Z}_n^k$ the Krawtchouk polynomial
$\mathcal{K}_k^{n;q}$ satisfies:
$$\mathcal{K}_k^{n;q}(|x|) = \sum_{y\in \mathbb{Z}_n^k: |y|=k}
\chi_y(x)~,$$ the eigenvalue corresponding to $y \in \mathbb{Z}_n^k$
is:
$$\mu(y) = \sum_{x\in S}\chi_x(y) =
\mathop{\sum_{t=0}^k}_{t\not\equiv
0~(\mathrm{mod}~p)}\sum_{x:|x|=t}\chi_x(y) =
\mathop{\sum_{t=0}^k}_{t\not\equiv 0~(\mathrm{mod}~p)}
\mathcal{K}_t^{k;q}(|y|)~.$$
\end{remark}
\begin{remark}
The upper bound on $\xal$ was derived from an asymptotic analysis of
the smallest eigenvalue $\mu_{n^k}$ of $G^k$. Tight results on
$\alpha(G^k)$ may be obtained by a careful analysis of the
expression in \eqref{eq-mu}. To illustrate this, we consider the
case $G=K_3$ and $p=3$. Combining the previous remark with
\eqref{eq-Kk-non-0-mod-p-sum} and \eqref{eq-xi-s-value}, we obtain
that the eigenvalues of $K_3^{k_{(3)}}$ are:
\begin{eqnarray} \mu_1 &=& \frac{2}{3}3^k -
\frac{2}{3}\sqrt{3}^k \cos(\frac{\pi
k}{2})~,\nonumber \\
\mu(s) &=& -\frac{2}{3} \sqrt{3}^k \cos(\frac{\pi k}{2} - \frac{2
\pi s}{3})~\mbox{ for }0 < s \leq k~.\end{eqnarray} Noticing that
$\mu(s)$ depends only on the values of $s\pmod{3}$ and $k \pmod{4}$,
we can determine the minimal eigenvalue of $G^k$ for each given
power $k$, and deduce that: $$\begin{array}{ll}
\displaystyle{\alpha(G^k) \leq
3^{k/2}}& \mbox{if }k \equiv 0\pmod{4}\\
\displaystyle{\alpha(G^k) \leq \frac{3^{k+1}}{3+2\cdot 3^{(k+1)/2}}
< \frac{1}{2}
3^{(k+1)/2}}& \mbox{if }k \equiv 1\pmod{2} \\
\displaystyle{\alpha(G^k) \leq \frac{3^k}{3+2\cdot 3^{k/2}} <
\frac{1}{2} 3^{k/2}}& \mbox{if }k \equiv 2\pmod{4}
\end{array}~,$$
matching the results obtained by the Delsarte linear programming
bound.

\end{remark}
\section{Ramsey subgraphs in large $p$-powers of any graph}\label{sec::ramsey}
In order to prove a poly-logarithmic upper bound on the clique sizes
of $p$-powers of a graph $G$, we use an algebraic argument, similar
to the method of representation by polynomials described in the
Section \ref{sec::general-bounds}. We note that the same approach
provides an upper bound on the size of independent sets. However,
for this latter bound, we require another property, which relates
the problem to strong graph products and to the Shannon capacity of
a graph.

The $k$-th {\em strong} power of a graph $G$ (also known as the {\em
and} power), denoted by $G^{\wedge k}$, is the graph whose vertex
set is $V(G)^k$, where two distinct $k$-tuples $u\neq v$ are
adjacent iff each of their coordinates is either equal or adjacent
in $G$:
$$ (u_1,\ldots,u_k) (v_1,\ldots,v_k) \in E(G^{\wedge k})~\mbox{
iff for all }i=1,\ldots,k:~ u_i = v_i \mbox{ or } u_i v_i \in
E(G)~.$$ In 1956, Shannon \cite{Shannon} related the independence
numbers of strong powers of a fixed graph $G$ to the effective
alphabet size in a zero-error transmission over a noisy channel.
Shannon showed that the limit of $\alpha(G^{\wedge
k})^{\frac{1}{k}}$ as $k\to\infty$ exists and equals $\sup_k
\alpha(G^{\wedge k})^{\frac{1}{k}}$, by super-multiplicativity; this
limit is denoted by $c(G)$, the Shannon capacity of $G$. It follows
that $c(G) \geq \alpha(G)$, and in fact equality holds for all
perfect graphs. However, for non-perfect graphs, $c(G)$ may exceed
$\alpha(G)$, and the smallest (and most famous) example of such a
graph is $C_5$, the cycle on $5$ vertices, where $\alpha(C_5)=2$ and
yet $c(C_5) \geq \alpha(C_5^{\wedge 2})^{\frac{1}{2}} = \sqrt{5}$.
The seemingly simple question of determining the value of $c(C_5)$
was solved only in 1979 by Lov\'{a}sz \cite{LovaszTheta}, who
introduced the $\vartheta$-function to show that $c(C_5)=\sqrt{5}$.

The next theorem states the bound on the clique numbers of
$G^{k_{(p)}}$, and relates the Shannon capacity of $\overline{G}$,
the complement of $G$, to bounds on independent sets of
$G^{k_{(p)}}$.
\begin{theorem}\label{thm-omega-alpha-bounds}
Let $G$ denote a graph on $n$ vertices and let $p \geq 2$ be a
prime. The clique number of $G^{k_{(p)}}$ satisfies:
\begin{equation}
  \label{eq-omega-g-k-bound}
 \omega(G^{k_{(p)}}) \leq \binom{k n + p-1}{p-1}~,
\end{equation}
and if $I$ is an independent set of both $G^{k_{(p)}}$ and
${\overline{G}}^{\wedge k}$, then:
\begin{equation}
  \label{eq-alpha-or-clique-bound}
  |I| \leq \binom{k n + \lfloor
\frac{k}{p}\rfloor}{\lfloor \frac{k}{p}\rfloor}~.
\end{equation}
Moreover, if in addition $G$ is regular then:
\begin{equation}\label{eq-omega-alpha-regular-bound}
 \omega(G^{k_{(p)}}) \leq \binom{k (n-1) + p}{p-1}~,~
  |I| \leq \binom{k(n-1) + \lfloor
\frac{k}{p}\rfloor+1}{\lfloor \frac{k}{p}\rfloor}~.\end{equation}
\end{theorem}
The above theorem implies that if $S$ is an independent set of
${\overline{G}}^{\wedge k}$, then any independent set $I$ of
$G^{k_{(p)}}[S]$, the induced subgraph of $G^{k_{(p)}}$ on $S$,
satisfies inequality \eqref{eq-alpha-or-clique-bound}. For large
values of $k$, by definition there exists such a set $S$ of size
roughly $c(\overline{G})^k$. Hence, there are induced subgraphs of
$G^{k_{(p)}}$ of size tending to $c(\overline{G})^k$, whose clique
number and independence number are bounded by the expressions in
\eqref{eq-omega-g-k-bound} and \eqref{eq-alpha-or-clique-bound}
respectively.

In the special case $G=K_n$, the graph ${\overline{G}}^{\wedge k}$
is an edgeless graph for any $k$, and hence:
$$\alpha(K_n^{k_{(p)}}) \leq \binom{k(n-1) +\lfloor \frac{k}{p}\rfloor+1}{\lfloor
\frac{k}{p}\rfloor} \leq
\left(\mathrm{e}p(n-1)+\mathrm{e}+o(1)\right)^{k/p}~,$$ where the
$o(1)$-term tends to $0$ as $k\to\infty$. This implies an upper
bound on $x_\alpha^{(p)}(K_n)$ which nearly matches the upper bound
of Theorem \ref{thm-xal-kn-bounds} for large values of $p$.
\begin{proof}
Let $g_1:V(G)\to \mathbb{Z}_p^m$ and $g_2:V(G)\to \mathbb{C}^m$, for
some integer $m$, denote two representations of $G$ by
$m$-dimensional vectors, satisfying the following for any (not
necessarily distinct) $u,v \in V(G)$:
\begin{equation}
  \label{eq-vector-representation}
  \left\{\begin{array}{ll}
    g_i(u)\cdot g_i(v) = 0 & \mbox{if }u v \in E(G)\\
    g_i(u)\cdot g_i(v) = 1 & \mbox{otherwise }
  \end{array}\right.~(i=1,2)~.
\end{equation}
It is not difficult to see that such representations exist for any
graph $G$. For instance, the standard basis of $n$-dimensional
vectors is such a representation for $G=K_n$. In the general case,
it is possible to construct such vectors inductively, in a way
similar to a Gram-Schmidt orthogonalization process. To see this,
define the lower diagonal $|V(G)|\times |V(G)|$ matrix $M$ as
follows:
$$ M_{k,i} = \left\{\begin{array}{cl}
-\sum_{j=1}^{i-1}M_{k,j} M_{i,j} & i < k,~v_i v_k \in E(G) \\
1-\sum_{j=1}^{i-1}M_{k,j} M_{i,j} & i < k,~v_i v_k \notin E(G) \\
1 & i = k \\
0 & i > k
\end{array}\right.~.$$
The rows of $M$ satisfy \eqref{eq-vector-representation} for any
distinct $u,v \in V(G)$, and it remains to modify the inner product
of any vector with itself into $1$ without changing the inner
products of distinct vectors. This is clearly possible over
$\mathbb{Z}_p$ and $\mathbb{C}$ using additional coordinates.

Consider $G^{k_{(p)}}$, and define the vectors $ w_u =
g_1(u_1)\circ\ldots\circ g_1(u_k)$ for $u=(u_1,\ldots,u_k)\in
V(G^k)$, where $\circ$ denotes vector concatenation. By definition:
$$w_u \cdot w_v \equiv k - |\{i:u_i v_i \in E(G)\}| \pmod{p}$$ for
any $u,v \in V(G^k)$, and hence, if $S$ is a maximum clique of
$G^k$, then $w_u \cdot w_v \not\equiv k \pmod{p}$ for any $u,v \in
S$. It follows that if $B$ is the matrix whose columns are $w_u$ for
$u\in S$, then $C=B^t B$ has values which are $k\pmod{p}$ on its
diagonal and entries which are not congruent to $k$ modulo $p$
anywhere else. Clearly, $\mathrm{rank}(C) \leq \mathrm{rank}(B)$,
and we claim that $\mathrm{rank}(B) \leq k n$, and that furthermore,
if $G$ is regular then $\mathrm{rank}(B) \leq k(n-1)+1$. To see
this, notice that, as the dimension of $\mathrm{Span}(\{g_1(u) : u
\in V\})$ is at most $n$, the dimension of the span of $\{w_u : u\in
G^k\}$ is at most $k n$. If in addition $G$ is regular, define $z =
\sum_{u\in V}g_1(u)$ (assuming without loss of generality that
$z\neq 0$), and observe that by \eqref{eq-vector-representation},
each of the vectors $w_u$ is orthogonal to the following $k-1$
linearly independent vectors:
\begin{equation}\label{eq-orth-vectors}\{z\circ(-z)\circ
\underline{0}^{\circ (k-2)},~\underline{0}\circ z\circ(-z)\circ
\underline{0}^{\circ(k-3)},\ldots,~ \underline{0}^{\circ(k-2)}\circ
z\circ(-z)\}~.\end{equation} Similarly, the vectors $ w'_u =
g_2(u_1)\circ\ldots\circ g_2(u_k)$ satisfy the following for any
$u,v \in V(G^k)$: $$w'_u \cdot w'_v = k - |\{i:u_i v_i \in
E(G)\}|~.$$ Let $I$ denote an independent set of $G^{k_{(p)}}$,
which is also an independent set of ${\overline{G}}^{\wedge k}$. By
the definition of ${\overline{G}}^{\wedge k}$, every $u,v\in I$
share a coordinate $i$ such that $u_i v_i \in E(G)$, and combining
this with the definition of $G^{k_{(p)}}$, we obtain:
$$ 0 < |\{i: u_i v_i \in E(G)\}| \equiv 0\pmod{p}~\mbox{ for any }u,v\in I~.$$
Therefore, for any $u \neq v\in I$: $$w'_u \cdot w'_v = k - t p
~\mbox{ for some }~t \in \{1,\ldots,\lfloor\frac{k}{p}\rfloor\}~,$$
and if $B'$ is the matrix whose columns are $w'_u$ for $u\in I$,
then $C'=B'^t B'$ has the entries $k$ on its diagonal and entries of
the form $k - tp$, $0<t\leq\lfloor\frac{k}{p}\rfloor$, anywhere
else. Again, the definition of $g_2$ implies that $\mathrm{rank}(C')
\leq kn$, and in case $G$ is regular, $\mathrm{rank}(C') \leq
k(n-1)+1$ (each $w'_u$ is orthogonal to the vectors of
\eqref{eq-orth-vectors} for $z=\sum_{u\in V}g_2(u)$).

Define the following polynomials:
\begin{equation}\label{eq-matrix-poly} f_1(x) = \mathop{\prod_{j \in
\mathbb{Z}_p}}_{j \not\equiv k(\mathrm{mod}~p)}
(j-x)~~,~~f_2(x)=\prod_{t=1}^{\lfloor\frac{k}{p}\rfloor}(k-tp-x)~.
\end{equation}
By the discussion above, the matrices $D,D'$ obtained by applying
$f_1,f_2$ on each element of $C,C'$ respectively, are non-zero on
the diagonal and zero anywhere else, and in particular, are of full
rank: $\mathrm{rank} (D) = |S|$ and $\mathrm{rank} (D') = |I|$.
Recalling that the ranks of $C$ and $C'$ are at most $kn$, and at
most $k(n-1)+1$ if $G$ is regular, the proof is completed by the
following simple Lemma of \cite{NogaExtComb}:
\begin{lemma}[\cite{NogaExtComb}]
Let $B = (b_{i,j})$ be an $n$ by $n$ matrix of rank $d$, and let
$P(x)$ be an arbitrary polynomial of degree $k$. Then the rank of
the $n$ by $n$ matrix $(P(b_{i,j}))$ is at most $\binom{k+d}{k}$.
Moreover, if $P(x) = x^k$ then the rank of $(P(b_{i,j}))$ is at most
$\binom{k+d-1}{k}$.\end{lemma}\end{proof}

For large values of $k$, the upper bounds provided by the above
theorem are:
\begin{eqnarray} \omega(H)&\leq& \binom{(1+o(1))k n}{p}~,\nonumber\\
\alpha(H) &\leq& \binom{(1+o(1))k n}{k/p}\nonumber~.
\end{eqnarray}
This gives the following immediate corollary, which states that
large $p$-powers of any nontrivial graph $G$ contain a large induced
subgraph without large homogenous sets.
\begin{corollary} Let $G$ be some fixed nontrivial graph and fix a prime
$p$.
\begin{enumerate}
\item Let $S$ denote a maximum clique of $G$, and set $\lambda = \log
\omega(G)=\log\alpha(\overline{G})$. For any $k$, the induced
subgraph of $G^{k_{(p)}}$ on $S^k$, $H=G^{k_{(p)}}[S^k]$, is a graph
on $N=\exp(k\lambda )$ vertices which satisfies:
$$ \omega(H) = O(\log^p N)~,~\alpha(H) \leq N^{(1+o(1))\frac{\log(np)+1}{p\lambda}}~.
$$
\item The above formula holds when taking $\lambda = \frac{\log
\alpha(\overline{G}^{\wedge \ell})}{\ell}$ for some $\ell\geq 1$
dividing $k$, $S$ a maximum clique of $\overline{G}^{\wedge \ell}$,
and $H = G^{k_{(p)}}[S^{k/\ell}]$. In particular, for sufficiently
large values of $k$, $G^{k_{(p)}}$ has an induced subgraph $H$ on
$N=\exp\left((1-o(1))k\log c(\overline{G})\right)$ vertices
satisfying: $$\omega(H) = O(\log^p N)~,~\alpha(H) \leq
N^{(1+o(1))\frac{\log (np)+1}{p \log c(\overline{G})}}~.$$
\end{enumerate}
\end{corollary}

\begin{remark}
In the special case $G=K_n$, where $n,p$ are large and $k > p$, the
bound on $\omega(K_n^k)$ is $\binom{(1+o(1))kn}{p}$ whereas the
bound on $\alpha(K_n^k)$ is $\binom{(1+o(1))kn}{k/p}$. Hence, the
optimal mutual bound on these parameters is obtained at $k=p^2$.
Writing $H=K_n^k$, $N=n^k=n^{p^2}$ and $p=n^c$ for some $c>0$, we
get:
$$ p = \sqrt{\frac{(2c+o(1))\log N}{\log\log N}}~,$$
and:
$$ \max\{\omega(H),\alpha(H)\} \leq
\left((1+o(1))\mathrm{e}pn\right)^p = \exp\left(
\left(\frac{1+c}{\sqrt{2c}}+o(1)\right)\sqrt{\log N \log\log
N}\right)~.$$ The last expression is minimized for $c=1$, and thus
the best Ramsey construction in $p$-powers of $K_n$ is obtained at
$p=n$ and $k=p^2$, giving a graph $H$ on $N$ vertices with no
independence set or clique larger than
$\exp\left((1+o(1))\sqrt{2\log N \log\log N}\right)$ vertices. This
special case matches the bound of the $FW$ Ramsey construction, and
is in fact closely related to that construction, as we next
describe.

The graph $FW_N$, where $N=\binom{p^3}{p^2-1}$ for some prime $p$,
is defined as follows: its vertices are the $N$ possible choices of
$(p^2-1)$-element sets of $[p^3]$, and two vertices are adjacent iff
the intersection of their corresponding sets is congruent to $-1$
modulo $p$. Observe that the vertices of the graph $K_n^{k_{(p)}}$
for $n=p$ and $k=p^2$, as described above, can be viewed as
$k$-element subsets of $[kn]$, where the choice of elements is
restricted to precisely one element from each of the $k$ subsets
$\{(j-1)n+1,\ldots,j n\}$, $j\in[k]$ (the $j$-th subset corresponds
to the $j$-th coordinate of the $k$-tuple). In this formulation, the
intersection of two sets corresponds to the number of common
coordinates between the corresponding $k$-tuples. As $k=p^2\equiv
0\pmod{p}$, it follows that two vertices in $K_p^{p^2_{(p)}}$ are
adjacent iff the intersection of their corresponding sets is not
congruent to $0$ modulo $p$. Altogether, we obtain that
$K_p^{p^2_{(p)}}$ is an induced subgraph of a slight variant of
$FW_N$, where the differences are in the cardinality of the sets and
the criteria for adjacency.

 Another relation between the two constructions is
the following: one can identify the vertices of $K_2^{p^3_{(p)}}$
with all possible subsets of $[p^3]$, where two vertices are
adjacent iff the intersection of their corresponding sets is not
congruent to $0$ modulo $p$. In particular, $K_2^{p^3_{(p)}}$
contains all the $(p^2-1)$-element subsets of $[p^3]$, a variant of
$FW_N$ for the above value of $N$ (the difference lies in the
criteria for adjacency).

We note that the method of proving Theorem
\ref{thm-omega-alpha-bounds} can be applied to the graph $FW_N$,
giving yet another simple proof for the properties of this well
known construction.
\end{remark}


\noindent\textbf{Acknowledgement} The authors would like to thank
Simon Litsyn and Benny Sudakov for useful discussions.

\end{document}